\documentclass[10pt]{article}
\usepackage[all]{xy}
\usepackage{amsmath}
\usepackage{amsfonts}
\usepackage{amssymb}
\usepackage{amscd}
\usepackage{amsthm}
\usepackage{latexsym}
\usepackage{amsbsy}

\begin{document}

\title{\bf{Weyl's Law and Connes' Trace Theorem  for  Noncommutative Two Tori}}

\author {FARZAD FATHIZADEH AND MASOUD KHALKHALI}

\date{}

\maketitle

\begin{abstract}

We prove the analogue of Weyl's law  for a noncommutative Riemannian manifold, 
namely the noncommutative two torus $\mathbb{T}_\theta^2$ equipped with a 
general translation invariant conformal structure and a Weyl conformal factor. This is 
achieved by studying the asymptotic distribution of the eigenvalues of the perturbed 
Laplacian  on $\mathbb{T}_\theta^2$. We  also prove the analogue of Connes' trace 
theorem by showing that the Dixmier trace and a noncommutative residue coincide
on pseudodifferential operators of order $-2$ on $\mathbb{T}_\theta^2$.

\end{abstract}



\section{Introduction}

A  celebrated    theorem of Weyl states that one can compute  the volume of a compact Riemannian 
manifold from the asymptotic growth  of the eigenvalues of its Laplacian. This  result,  known as 
Weyl's  law,  has  been the starting point of numerous developments and  conjectures  in 
mathematics and physics. Some of these results are now presented in the context of spectral geometry, 
inverse problems  and quantum mechanics. It is pointed out in \cite{guil} that a subtle asymptotic identity 
seems to be the only intrinsic link between a classical observable $h$   and its  quantization $H$.  That is, the eigenvalue counting function  
for $H$ defined by   $N (\Lambda) =  \#  \{ \lambda  \leq  \Lambda \},$   and the phase space volume 
$\text{Vol}  (h\leq \Lambda)$ satisfy the asymptotic identity
\[
N(\Lambda) \sim \gamma \,\textnormal{Vol}(h \leq \Lambda) \qquad \textnormal{as} \qquad \Lambda \to \infty, 
\]
where $\gamma$ is a constant that does not depend on  $h$ or $H$. Here $h$ is  a  real-valued function 
 on a symplectic  phase space and  $H$ is a selfadjoint 
operator  on a Hilbert space.

Like similar results in spectral geometry which encode geometric information in terms of spectral data, Weyl's law  
has  played an important role in the development of metric aspects of noncommutative   geometry 
\cite{con-1, con0, con1, con1.5,  conmos1.5, chacon,  conmar}.  A noncommutative geometric space is described by 
a spectral triple $(\mathcal{A}, \mathcal{H}, D)$. Here $\mathcal{A}$ is an algebra acting by bounded operators on a 
Hilbert space $\mathcal{H}$, and $D$ is an unbounded 
selfadjoint operator acting in $\mathcal{H}$ such that the commutators $[D, a]$ are bounded 
for all $a \in \mathcal{A}$ \cite{con1, conmos1.5, con1.5, conmar}. 
Under suitable summability and regularity conditions,  the geometric invariants such as volume and curvature of 
a noncommutative space $(\mathcal{A}, \mathcal{H}, D)$ with metric dimension $d$ are obtained by considering 
small time asymptotic expansions of the form  
\[ 
\text{Trace} (a e^{-t  D^2}) \sim  \sum_{n=0}^\infty a_{n} ( a, D^2) t^{\frac{n-d}{2}}, \qquad  a \in \mathcal{A}.
\] 
There is an equivalent formulation of such geometric invariants in terms of spectral zeta functions. 

The main results of this paper are Theorem \ref{weyllaw} and Theorem \ref{connestrace} in which we prove  
the analogue of Weyl's law and Connes' trace theorem \cite{con0} for  the noncommutative two torus 
$\mathbb{T}_\theta^2$. In  Section \ref{prel} we recall the construction of a spectral triple on 
$\mathbb{T}_\theta^2$  via a complex 
number $\tau$ in the upper half-plane  and a  positive invertible 
element $k \in C^\infty(\mathbb{T}_\theta^2)$. The complex parameter $\tau$ defines a translation 
invariant conformal structure on $\mathbb{T}_\theta^2$ and the element $k$ plays the role of a Weyl conformal factor 
\cite{contre, fatkha}. In Section 
\ref{ncpseudos} we recall Connes'  pseudodifferential calculus \cite{con} for the canonical 
dynamical system defining the noncommutative two torus.  Using this 
pseudodifferential calculus and heat kernel techniques \cite{gil, contre, fatkha}, in Section \ref{weyllawsec} we  
find the first coefficient in the small time heat kernel expansion of  the perturbed Laplacian on 
$\mathbb{T}_\theta^2$.  The volume of  $\mathbb{T}_\theta^2$ with respect to the perturbed volume 
form manifests itself in this coefficient, and we use a Tauberian theorem to find its relation with the asymptotic distribution of the 
eigenvalues of the Laplacian. This establishes the analogue of Weyl's law for $\mathbb{T}_\theta^2$. In Section \ref{connestracesec} 
we prove the analogue of Connes' trace theorem \cite{con0} for 
$\mathbb{T}_\theta^2$. That is, we show that the Dixmier trace \cite{dix} and a noncommutative 
residue defined in \cite{fatwon} coincide on pseudodifferential operators of order $-2$ on 
$\mathbb{T}_\theta^2$.

Recently there has been much progress in understanding the local  differential geometry of the  noncommutative  
two torus equipped with a general metric. In \cite{fatkha}, the  Gauss-Bonnet theorem of \cite{contre} is  extended 
to  full generality following the pioneering work \cite{cohcon} in the subject. The much more delicate 
question of computing the scalar curvature of this noncommutative manifold is now also fully settled in 
\cite{conmos2,fatkha1}. In \cite{conmos2}, among  other results,  one can also find a  closed formula 
for the Ray-Singer analytic torsion in terms of the Dirichlet quadratic form,  and the 
corresponding evolution equation for the metric is shown to give  the appropriate analogue of Ricci curvature. 
A different approach to the question of Ricci curvature is proposed in \cite{bhumar}. 

The first author would like  to thank IHES for  kind support and excellent environment during his visit in Summer 
2011 where part of this work was carried out. We would like to thank the referees for carefully reading the paper 
and making several useful suggestions.

\section{Preliminaries} \label{prel}

Let $\Sigma$ be a closed, oriented,  2-dimensional smooth manifold equipped with a Riemannian metric $g$. 
The spectrum of the Laplacian $\triangle_{g} = d^*d$, where $d$ is the de Rham differential operator acting on 
smooth functions on $\Sigma$, encodes important geometric and topological information about  $(\Sigma, g)$. 
In fact there is an asymptotic expansion for the trace of the heat kernel $e^{-t \triangle_g}$ which is of the form 
\begin{equation} \label{claasym}
\text{Trace} (e^{-t \triangle_g}) \sim t^{-1} \sum_{n=0}^\infty a_{2n} (\triangle_g) t^n \qquad (t \to 0),
\end{equation} 
with $a_{0} (\triangle_g) =  \frac{1}{4 \pi}\text{Vol}(\Sigma) $ and $a_{2} (\triangle_g) = \frac{1}{24 \pi} K(\Sigma)$. Here  
$\text{Vol}(\Sigma)$ and $K(\Sigma)$ respectively denote the volume and total curvature of the Riemannian manifold.

On the other hand $\text{Trace}(e^{-t \triangle_g})$ depends only on the spectrum of the Laplacian since  
\begin{equation} \label{trofheatker}
\text{Trace}(e^{-t \triangle_g}) = \sum e^{-t \lambda_j}, 
\end{equation} 
where the summation is over all eigenvalues $\lambda_j$ of $\triangle_g$ counted with multiplicity. 
Comparing relations \eqref{claasym}, \eqref{trofheatker}, and using Karamata's  Tauberian 
theorem \cite{bgv} (page 91) the term $a_0 =  \frac{1}{4 \pi} \text{Vol}(\Sigma) $ in the above asymptotic 
expansion is seen to depend only on the asymptotic growth of the eigenvalues 
$\lambda_j$. More precisely, one obtains 
\[  
\lim_{j \to \infty} \frac{\lambda_j}{j} = \frac{4 \pi}{ \text{Vol}(\Sigma)}.
\] 
This provides a proof for the celebrated Weyl's law which states that one can hear the volume of 
a Riemannian manifold from the asymptotic distribution of the eigenvalues of its Laplacian. 

In the following subsections we recall the definition of the noncommutative two tours and the construction of 
a spectral triple on it. The spectral triple will encode the metric information corresponding to general translation invariant 
conformal structures on $\mathbb{T}_\theta^2$ and a Weyl conformal factor by means of which one can vary 
inside the conformal class \cite{cohcon, contre, fatkha}.

\subsection{The irrational rotation algebra}

Let  $\theta$ be an irrational number. Recall that the  irrational rotation $C^*$-algebra $A_{\theta}$ is, 
by definition, the universal unital $C^*$-algebra  generated by two unitaries $U, V$ satisfying 
\[
VU=e^{2 \pi i \theta} UV.
\] 
One usually thinks of $A_{\theta}$ as the algebra of continuous functions on the noncommutative 2-torus 
$\mathbb{T}_{\theta}^2$. There is a continuous action of the ordinary torus $\mathbb{T}^2$, 
$\mathbb{T}= \mathbb{R}/2\pi \mathbb{Z}$, on $A_{\theta}$ by $C^*$-algebra automorphisms  
$\{ \alpha_s\}_{s\in \mathbb{R}^2}$ defined by 
\[
\alpha_s(U^mV^n)=e^{is.(m,n)}U^mV^n, \qquad m, n \in \mathbb{Z}.
\] 
The space of smooth elements for this action, denoted by $A_\theta^\infty$, consists of  those elements 
$a \in A_{\theta}$ for which the map $s \mapsto \alpha_s (a)$  is a smooth map from $\mathbb{R}^2$ to $A_\theta$. 
It is a dense subalgebra of $A_{\theta}$  which can be alternatively described as the algebra of elements of the form 
\[
\sum_{m,n\in \mathbb{Z}}a_{m,n}U^mV^n
\] 
with rapidly decreasing coefficients, \emph{i.e.} for any non-negative integers $k,q,$ we have 
\[ 
\sup_{m,n\in \mathbb{Z}}
\big ( (1+|m|)^k (1+|n|)^q|a_{m,n}| \big ) < \infty.
\]
There is a unique normalized trace   $\mathfrak{t}$  on $A_{\theta}$ whose restriction to the smooth 
elements is given by 
\[
\mathfrak{t}\big (\sum_{m,n\in \mathbb{Z}}a_{m,n}U^mV^n \big )=a_{0,0}.
\]
This trace functional is positive and faithful, and is used to define an inner product on $A_\theta$:
\begin{equation} \label{innerp}
\langle a, b \rangle = \mathfrak{t}(b^*a), \qquad a,b \in A_{\theta}.  
\end{equation}
Completing $A_\theta$ with respect to this inner product, we obtain a Hilbert space which will be 
denoted by $\mathcal{H}_0$.
 
The infinitesimal generators of the above action of $\mathbb{T}^2$ on $A_{\theta}$ are the 
derivations  $\delta_1,  \delta_2: A_{\theta}^{\infty} \to A_{\theta}^{\infty}$ defined by 
\[
\delta_1(U^mV^n)=mU^mV^n,  \qquad  \delta_2(U^mV^n)=nU^mV^n, \qquad  m,n \in \mathbb{Z}.
\]
In fact, $\delta_1, \delta_2$ are analogues of the differential operators
$\frac{1}{i}\partial/\partial x, \frac{1}{i}\partial/\partial y$ acting on smooth functions on
the ordinary two torus. We have $\delta_j(a^*)= -\delta_j(a)^* $ for $j=1, 2,$ and all $a\in A_{\theta}^{\infty}$.
Moreover, since $\mathfrak{t} \circ \delta_j =0$ for $j=1, 2$, we have the
analogue of integration by parts:
\[ 
\mathfrak{t}(a\delta_j(b)) = -\mathfrak{t}(\delta_j(a)b), \qquad  a,b \in A_{\theta}^{\infty}. 
\]
As unbounded operators on $\mathcal{H}_0$,  the derivations $\delta_1, \delta_2$ are formally 
selfadjoint and have unique selfadjoint extensions.

\subsection{Conformal structures  on $A_\theta$} \label{conf}
In this subsection we recall the analogues of complex and conformal structures for noncommutative two tori. Our definitions are motivated 
by their commutative counterparts. To  any complex number $\tau$ in the upper half-plane, $\Im(\tau) > 0$, we can 
associate a complex structure on the noncommutative two torus by defining the operators
\[ 
\partial = \delta_1 + \bar \tau \delta_2, \qquad \partial^*=  \delta_1 + \tau \delta_2. 
\]
These derivations  on $A_{\theta}^{\infty}$  are analogues of Dolbeault operators on Riemann surfaces. 
We note that   $\partial$ is an unbounded densely defined  operator on   $\mathcal{H}_0$ and $\partial^*$ is
 its formal adjoint. 
 
It is explained in \cite{con1} (\S VI. 2) that the information on the conformal structure represented by $\tau$ 
is encoded in the positive Hochschild $2$-cocycle on $A_\theta^\infty$ given by
\[\psi(a, b, c) = - \mathfrak{t}(a \partial(b) \partial^*(c)), \qquad a, b, c \in A_\theta^\infty.\]
Therefore, the analogue of the space of $(1,0)-$forms on the ordinary two torus
is defined to be the Hilbert space completion of the space of finite sums $\sum a \partial b$, $a,b \in A_\theta^\infty,$ 
with respect to the inner product defined by
\[
\langle a \partial(b), a'\partial(b') \rangle = \mathfrak{t}((a')^*a \partial(b)\partial(b')^*), \qquad a, b, a', b' \in A_{\theta}^\infty.
\]
This Hilbert space will be denoted by $\mathcal{H}^{(1, 0)}$.

In order to vary inside the conformal class represented by $\tau$, we choose a 
selfadjoint element $h \in A_{\theta}^{\infty}$ (\emph{cf.} \cite{contre}), and define a positive linear functional
$\varphi$ on $A_{\theta}$ by 
\[
\varphi(a)= \mathfrak{t}(ae^{-h}), \qquad a \in A_{\theta}.
\]
This linear functional is not a trace. However,
it is a twisted trace and satisfies the KMS condition at $\beta = 1$ for the
1-parameter group $\{\sigma_t \}_{t \in \mathbb{R}}$ of inner automorphisms
\[
\sigma_t= \Delta^{-it},
\] 
where the modular operator $\Delta$  is defined  by
(\emph{cf.} \cite{contre})
\[
\Delta(a)=e^{-h}ae^{h}, \qquad a \in A_\theta.
\]
We define a new inner product $\langle \, , \, \rangle_{\varphi}$ on $A_{\theta}$ by
\[ 
\langle a,b  \rangle_{\varphi} = \varphi(b^*a), \qquad a,b \in A_{\theta}. 
\]
The Hilbert space obtained from completing $A_{\theta}$ with respect to this 
inner product will be denoted by $\mathcal{H}_{\varphi}$.

\subsection{A spectral triple on $A_{\theta}$ } \label{spectripsubsec}

Following the above notations, let $\partial_{\varphi}$ be the same operator as  $\partial$  
but viewed as an unbounded operator from $\mathcal{H}_{\varphi}$ to $\mathcal{H}^{(1,0)}$. 
Consider the left action of $A_{\theta}$ on the 
Hilbert space
\[   \mathcal{H} = \mathcal{H}_\varphi \oplus \mathcal{H}^{(1,0)}\]
defined by $a (b, c)=(ab, ac)$ for all $a\in A_{\theta}, b \in  \mathcal{H}_{\varphi},$ and $c \in \mathcal{H}^{(1,0)}$. 
The operator
\[ D=
\left(\begin{array}{c c}
0 & \partial_\varphi^* \\
\partial_\varphi & 0
\end{array}\right)
: \mathcal{H} \to \mathcal{H}\]
is an unbounded selfadjoint operator and  one checks that this data defines a spectral triple  which we denote by 
$(A_\theta, \mathcal{H}, D)$  \cite{contre}.
This is an even spectral triple with the grading operator 
given by 
\[
\gamma =  \left(\begin{array}{c c}
 1 & 0  \\
0 & -1
\end{array}\right) : \mathcal{H} \to \mathcal{H}.
\] 

The above spectral triple endows the noncommutative two torus with a metric structure in the sense of 
noncommutative geometry (cf. \cite{con1}, \S VI). More precisely it  encodes the metric information 
represented by the complex structure and the Weyl conformal factor introduced in Subsection \ref{conf}. 
Apart from this paper,  this spectral triple  played an important role in recent work on the Gauss-Bonnet 
theorem and scalar curvature for noncommutative two tori \cite{contre, fatkha, conmos2, fatkha1}. 

The  full Laplacian \begin{equation} 
D^2 =  \left(\begin{array}{c c}
\partial_\varphi^*\partial_\varphi & 0  \\
0 & \partial_\varphi \partial_\varphi^*
\end{array}\right) \nonumber
\end{equation}
consists of the Laplacian on $\mathcal{H}_\varphi$ and the Laplacian on $\mathcal{H}^{(1, 0)}$. 
Since these operators have the same spectrum  away from $0$, for the purpose of this paper,  we only 
work with the first Laplacian 
\[
\triangle':=\partial_\varphi^* \partial_\varphi: \mathcal{H}_\varphi \to \mathcal{H}_\varphi. 
\]
In the flat case, when $h=0$,  one can easily see that this Laplacian reduces to the flat Laplacian 
\[
\triangle:=\partial^* \partial = \delta_1^2 + 2 \Re(\tau) \delta_1\delta_2 +
|\tau|^2\delta_2^2. 
\]

In general, $ \triangle':\mathcal{H}_\varphi \to \mathcal{H}_\varphi$ 
is a  complicated operator. However, it can be shown that  it is anti-unitarily equivalent to the operator
\[k \triangle k :\mathcal{H}_0 \to \mathcal{H}_0,\] 
where $k:= e^{h/2}$ acts on $\mathcal{H}_0$ by left multiplication ({\it cf.} \cite{contre, fatkha}). In particular 
these two operators have the same spectrum. This observation  plays a crucial role  in this paper.

\section{Pseudodifferential Operators on $\mathbb{T}_\theta^2$} \label{ncpseudos}

Connes' pseudodifferential calculus \cite{con} may be used to compute spectral invariants of 
noncommutative tori as in  \cite{contre,conmos2,fatkha,fatkha1}. Our proof of Weyl's law and Connes' trace 
theorem for noncommutative tori is also based on this calculus. 

A differential operator of order $n\geq 0$  on the noncommutative two torus is an operator of the 
form
\[
\sum_{ j_1, j_2} a_{j_1, j_2} \delta_1^{j_1} \delta_2^{j_2}: A_\theta^\infty \to A_\theta^\infty, 
\]
where $a_{j_1, j_2} \in A_\theta^\infty$ and $j_1, j_2$ are non-negative integers such that $j_1+j_2 \leq n$. 
The elements $a_{j_1, j_2}$ act by left multiplication. Using operator-valued symbols, the above  
notion of differential operators can be generalized to the notion of pseudodifferential operators.  Here 
we briefly recall these symbols and the corresponding pseudodifferential calculus. We shall use the notation 
$\partial_1=\frac{\partial}{\partial \xi_1}$, $\partial_2=\frac{\partial}{\partial \xi_2}$ for 
partial differentiation with respect to the coordinates of $\mathbb{R}^2$. 

\newtheorem{main1}{Definition}[section]
\begin{main1}
Let $n$ be an integer. A  smooth map $\rho: \mathbb{R}^2 \to A_{\theta}^{\infty}$ is  said
to be a symbol of order $n$  if for all non-negative integers $i_1, i_2, j_1,
j_2 $ we have
\[ 
||\delta_1^{i_1} \delta_2^{i_2} \partial_1^{j_1} \partial_2^{j_2} \rho(\xi) || \leq c (1+|\xi|)^{n-j_1-j_2},
\]
where $c$ is a constant, and if there exists a smooth map $k: \mathbb{R}^2 \setminus \{ 0\} \to
A_{\theta}^{\infty}$ such that
\[
\lim_{\lambda \to \infty} \lambda^{-n} \rho(\lambda\xi_1, \lambda\xi_2) = k (\xi_1, \xi_2), \qquad (\xi_1, \xi_2) \in \mathbb{R}^2 \setminus \{ 0\}.
\]
Here $||\cdot||$ denotes the $C^*$-norm of $A_\theta$. The space of symbols of order $n$ is denoted by $S_n$.
\end{main1}
To a symbol $\rho \in S_n$  one can associate an operator $P_{\rho}$ on $A_{\theta}^{\infty}$ defined by 
\begin{equation} \label{pseudofourierdef}
P_{\rho}(a) = (2 \pi)^{-2} \int_{\mathbb{R}^2} \int_{\mathbb{R}^2} e^{-is \cdot \xi} \rho(\xi) \alpha_s(a) \,ds \, d\xi, \qquad a \in A_\theta^\infty.
\end{equation}
The operator $P_{\rho}$ is said to be a pseudodifferential operator of order $n$. For
example, the differential operator $\sum_{j_1+j_2 \leq n } a_{j_1,j_2}
\delta_1^{j_1} \delta_2^{j_2}$ is associated with the symbol $\sum_{j_1+j_2 \leq n } a_{j_1,j_2}
\xi_1^{j_1} \xi_2^{j_2}$ via the above formula.

\newtheorem{symbolequivalence}[main1]{Definition}

\begin{symbolequivalence}
Two symbols $\rho$, $\rho'\in S_k$  are said to be equivalent if  $\rho-\rho'$ is in
$S_n$ for all integers $n$. The equivalence of  symbols will be denoted by  $\rho \sim \rho'$.
\end{symbolequivalence}

The following  Proposition in  \cite{con} shows that the space of pseudodifferential operators on $\mathbb{T}^2_\theta$
is an algebra and gives an explicit formula for  the symbol of the product of two 
pseudodifferential operators up to the above equivalence relation. It also shows that 
the formal  adjoint of a pseudodifferential operator, with respect to the inner
product defined on $A_{\theta}^{\infty} $ by \eqref{innerp}, is again a pseudodifferential
operator and gives a formula for its symbol.

\newtheorem{symbolcalculus}[main1]{Proposition}

\begin{symbolcalculus} \label{symbolcalculus}
Let $P$ and $Q$ be pseudodifferential operators with the symbols
$\rho$ and $\rho'$ respectively. Then the adjoint $P^*$ and
the product $PQ$ are pseudodifferential operators with the following
symbols
\[
\sigma(P^*) \sim \sum_{\ell_1, \ell_2 \geq 0} \frac{1}{\ell_1! \ell_2!}
\partial_1^{\ell_1} \partial_2^{\ell_2}\delta_1^{\ell_1}\delta_2^{\ell_2}
(\rho(\xi))^*,
\]
\[
\sigma (P Q) \sim \sum_{\ell_1, \ell_2 \geq 0} \frac{1}{\ell_1! \ell_2!}
\partial_1^{\ell_1} \partial_2^{\ell_2}(\rho (\xi))
\delta_1^{\ell_1}\delta_2^{\ell_2} (\rho'(\xi)).
\]
\end{symbolcalculus}

We also recall the notion of ellipticity for these pseudodifferential operators:

\newtheorem{elliptic}[main1]{Definition}

\begin{elliptic} \label{elliptic}
Let $\rho$ be a symbol of order $n$. It is said to be elliptic if $\rho(\xi)$ is invertible for $\xi \neq 0$,
and if there exists a constant $c$ such that
\[ || \rho(\xi)^{-1} || \leq c (1+|\xi|)^{-n} \]
for sufficiently large $|\xi|.$
\end{elliptic}
The flat Laplacian $\triangle=\delta_1^2 + 2 \Re(\tau) \delta_1\delta_2 +
|\tau|^2\delta_2^2$ is an example of an  elliptic operator  on $\mathbb{T}_\theta^2$.

\section{Weyl's Law for $\mathbb{T}_\theta^2$} \label{weyllawsec}

The resolvent of the Laplacian of a  closed Riemannian manifold  can be approximated  by 
pseudodifferential operators. Combining this with the Cauchy integral formula one can find a 
small time asymptotic expansion for the trace of the heat kernel of the Laplacian \cite{gil}. As 
stated in \cite{contre} ({\it cf.} also \cite{fatkha, conmos2, fatkha1}), one can use Connes' pseudodifferential calculus 
to find similar asymptotic expansions in noncommutative settings. We recall the main idea of this technique and will use it to establish 
the analogue of Weyl's law for $\mathbb{T}_\theta^2$. 

Let $\triangle'$ be the perturbed Laplacian on $A_\theta$ introduced in Subsection \ref{spectripsubsec}. 
Using the Cauchy integral formula, one has
\[ 
e^{-t\triangle'} = \frac{1}{2\pi i} \int_C e^{-t \lambda} (\triangle' - \lambda)^{-1} \, d \lambda, 
\]
where $C$ is a curve in the complex plane that goes around the non-negative
real axis  in such a way that
\[ 
e^{-t s} = \frac{1}{2\pi i} \int_C e^{-t \lambda} (s - \lambda)^{-1} \, d \lambda, \qquad s \geq 0.
\]
Appealing to this formula, one can use Connes' pseudodifferential calculus to  employ similar
arguments to those in \cite{gil} and derive an asymptotic expansion of the form 
\begin{equation} \label{asympl}
\text{Trace}(e^{-t \triangle'}) \sim t^{-1} \sum_{n=0}^{\infty} B_{2n} (\triangle')
t^n \qquad (t \to 0). \nonumber
\end{equation} 
That is, one can approximate $(\triangle' - \lambda)^{-1}$ 
by pseudodifferential operators $B_\lambda$ whose symbols are of the form 
\[
b_0(\xi, \lambda) +  b_1(\xi, \lambda) +  b_2(\xi, \lambda) + \cdots, 
\]
where for $j=0, 1, 2, \dots$, $b_j(\xi, \lambda)$ is a symbol of order $-2-j$. 

Since $\text{Trace}(e^{-t \triangle'})$  depends only on the eigenvalues of $\triangle'$, for the purpose of deriving the above asymptotic expansion one can work with 
the operator $k \triangle k$ which as shown in \cite{contre, fatkha} is anti-unitarily equivalent to $\triangle'$. 
The symbol of $k \triangle k$ as a pseudodifferential operator \cite{fatkha} is equal to
$a_2(\xi)+a_1(\xi)+a_0(\xi)$, where
\[ 
a_2(\xi)= \xi_1^2k^2+|\tau|^2\xi_2^2k^2+2 \Re(\tau) \xi_1\xi_2k^2,  
\]
\[
a_1(\xi)=2\xi_1k\delta_1(k) + 2|\tau|^2\xi_2k\delta_2(k) +
2 \Re(\tau)\xi_1k\delta_2(k)+2 \Re(\tau)\xi_2k\delta_1(k),
\]
\[
a_0(\xi)= k\delta_1^2(k)+ |\tau|^2k\delta_2^2(k) + 2 \Re(\tau)k\delta_1\delta_2(k).
\]
Therefore, using the calculus of symbols described in Section \ref{ncpseudos},  one has to solve 
the equation 
\begin{eqnarray} 
(b_0+b_1+b_2+\cdots) ((a_2-\lambda)+a_1+a_0) \sim 1. \nonumber
\end{eqnarray}
Here, $\lambda$ is treated as a symbol of order $2$ and we let $a_2'=a_2-\lambda, a_1'=a_1, a_0'=a_0$. 
Then the above equation yields
\[ 
\sum_{\substack{j, \ell_1, \ell_2 \geq 0,\\ k=0, 1, 2}} \frac{1}{\ell_1! \ell_2!}
\partial_1^{\ell_1} \partial_2^{\ell_2}(b_j)
\delta_1^{\ell_1}\delta_2^{\ell_2} (a_k') \sim 1.
\]
Comparing symbols of the same order on both sides, one concludes that 
\begin{equation} 
b_0=a_2'^{-1}=(a_2 - \lambda)^{-1}=(\xi_1^2k^2+|\tau|^2\xi_2^2k^2+
2\Re(\tau)\xi_1\xi_2k^2 - \lambda)^{-1}, \nonumber
\end{equation}
and
\[b_n  = - \sum_{\substack{2+j+\ell_1+\ell_2-k=n, \\ 0 \leq j <n, \, 0 \leq k \leq 2}} \frac{1}
{\ell_1! \ell_2!} \partial_1^{\ell_1} \partial_2^{\ell_2}(b_j) \delta_1^{\ell_1}\delta_2^{\ell_2} (a_k) b_0, \qquad n>0.\]
Similar to \cite{gil} one can use these symbols to approximate $e^{-t k \triangle k}$ with suitable infinitely 
smoothing operators and derive the desired asymptotic expansion. We record this result, already stated in \cite{contre}, 
in the following proposition.

\newtheorem{asympexpa}[main1]{Proposition} 

\begin{asympexpa} \label{asympexpa}
Let $\triangle'$ be the perturbed Laplacian on $A_\theta$ as above. There is an asymptotic expansion 
\begin{equation} \label{expansion}
\textnormal{Trace}(e^{-t \triangle'}) \sim t^{-1} \sum_{n=0}^{\infty} B_{2n} (\triangle') t^n \qquad (t \to 0),
\end{equation}
where  for each $n=0, 1, 2, \dots,$ 
\[B_{2n}(\triangle') = \frac{1}{2 \pi i} \int \int_C e^{-\lambda}  \mathfrak{t}(b_{2n}(\xi, \lambda)) \, d \lambda\, d \xi.\]
\end{asympexpa}

The following theorem is the analogue of Weyl's law for the noncommutative two torus equipped with a
general metric. In fact we show that the first term in the above asymptotic expansion is intimately
related to a natural definition for the volume of $(\mathbb{T}_\theta^2, \tau, k)$, and  by means of a Tauberian
theorem the relation between the volume and  the eigenvalue counting function of $\triangle'$
 is derived. We recall that $\tau$ is a complex parameter in the upper 
half-plane associated to which is a conformal structure on $\mathbb{T}_\theta^2$, and $k=e^{h/2}$ is an invertible 
selfadjoint element of $A_\theta^\infty$ that plays the role of a Weyl conformal factor.

\newtheorem{weyllaw}[main1]{Theorem} 

\begin{weyllaw} \label{weyllaw}
For any positive real number $\lambda$, let $N(\lambda)$ denote the number of eigenvalues of the perturbed Laplacian
$\triangle'$ that are less than $\lambda$. Then as $\lambda \to \infty$, we have
\begin{equation} \label{asympcount}
N(\lambda) \sim  \frac{\pi}{\Im(\tau)} \varphi(1) \lambda.
\end{equation}

\begin{proof}

Using Proposition \ref{asympexpa}  and passing to the coordinates \cite{fatkha}
\[
\xi_1= r \cos \theta - r \frac{\Re(\tau)}{\Im(\tau)} \sin \theta, \qquad \xi_2=\frac{r}{\Im(\tau)}\sin \theta ,
\]
where $\theta$ ranges from 0 to $2\pi$ and $r$ goes from 0 to $\infty$, we can compute  the first coefficient
in the asymptotic expansion \eqref{expansion} for the trace of the heat kernel of $\triangle'$ directly:
\[ 
B_0(\triangle') = \frac{1}{2 \pi i} \int \int_C e^{-\lambda}  \mathfrak{t}(b_0(\xi, \lambda)) \, d \lambda \, d \xi=
\frac{\pi}{\Im(\tau)}  \mathfrak{t}(k^{-2}). 
\]
Note that 
\[b_0=a_2'^{-1}=(\xi_1^2k^2+|\tau|^2\xi_2^2k^2+ 2\Re(\tau)\xi_1\xi_2k^2 - \lambda)^{-1}.\]
Now the asymptotic behavior of the eigenvalue counting function $N$ can be determined as follows.
It follows from \eqref{expansion} that
\[ 
\lim_{t \to 0^+} t \sum e^{-t \lambda_j} = B_0(\triangle') = \frac{\pi}{\Im(\tau)} \mathfrak{t}(k^{-2}),
\]
where the summation is over all eigenvalues $\lambda_j$ of $\triangle'$. Then it follows immediately from
 Karamata's Tauberian theorem \cite{bgv} (page 91) that the eigenvalue counting function $N$ satisfies
\[ 
N({\lambda}) \sim \frac{\pi}{\Im(\tau) \Gamma(2)} \mathfrak{t}(k^{-2}) \lambda=  \frac{\pi}{\Im(\tau)} 
\mathfrak{t}(k^{-2}) \lambda \qquad \textnormal{as} \qquad \lambda \to \infty.
\]

\end{proof}

\end{weyllaw}

A corollary of this theorem is that $(1+\triangle')^{-1}$ is in the domain of the Dixmier trace. 
Before stating the corollary we quickly review the Dixmier trace and the  noncommutative integral, 
following  \cite{con1}.  

We denote the ideal of compact operators on a  Hilbert space $\mathcal{H}$ by $\mathcal{K}(\mathcal{H})$.
For  any $T \in \mathcal{K}(\mathcal{H})$,  let $\mu_n(T), n=1, 2, \dots,$ denote the sequence of eigenvalues of 
its absolute value $|T|=(T^*T)^{\frac{1}{2}}$ written in  decreasing order with multiplicity:
\[
\mu_1(T) \geq \mu_2(T) \geq \cdots  \geq 0.
\]
The Dixmier trace is a trace functional on an ideal of compact operators $\mathcal{L}^{1,\infty } (\mathcal{H})$ defined as
\[
\mathcal{L}^{1,\infty } (\mathcal{H})= \big \{T \in \mathcal{K}(\mathcal{H}); \qquad \sum_{n=1}^N \mu_n (T)=O \,(\text{log} N ) \big \}.
\]
This ideal of operators is equipped with a natural norm:
\[
||T||_{1, \infty} := \sup_{N \geq 2} \frac{1}{\log N}\sum_{n=1}^N \mu_n(T), \qquad T \in \mathcal{L}^{1,\infty } (\mathcal{H}).
\]
It is clear that trace class operators are automatically in $\mathcal{L}^{1,\infty} (\mathcal{H})$.
The Dixmier trace of an operator $T \in \mathcal{L}^{1,\infty}(\mathcal{H})$ measures the  logarithmic divergence of
its ordinary trace. More precisely,  for any positive operator $T\in \mathcal{L}^{1,\infty } (\mathcal{H})$ we are interested 
in  the limiting behavior  of the  sequence
\[
\frac{1}{\log N}\sum_{n=1}^N \mu_n(T), \qquad N = 2, 3, \dots \,.
\]
While by our assumption this sequence is bounded,  its usual limit may not exist and 
must be replaced by a suitable generalized limit. The limiting procedure is carried out by 
means of a state on a $C^*$-algebra. Recall that a state on a  $C^*$-algebra is a non-zero positive linear 
functional on the algebra.

To define the Dixmier trace of a positive operator $T \in \mathcal{L}^{1,\infty}(\mathcal{H})$, consider  
the partial trace
\[
\text{Trace}_{N} (T) =\sum_{n=1}^N \mu_n (T), \qquad N=1, 2, \dots,
\]
and its piecewise affine interpolation denoted by $ \text{Trace}_{r} (T)$ for  $r \in [1, \infty).$ 
Then let 
\[
\tau_{\Lambda} (T): =  \frac{1}{\log  \Lambda} \int_e^{\Lambda} \frac{ \text{Trace}_{r} (T)}{\log r}\frac{dr}{r}, \qquad 
\Lambda \in [e, \infty), 
\]
be the  Ces\`{a}ro mean of the function $\text{Trace}_{r} (T)/\log r$ over the multiplicative group $\mathbb{R}^*_+$.
Now choosing a normalized state $\omega : C_b[e, \infty) \to \mathbb{C}$  on the algebra of bounded continuous functions
on  $[e, \infty)$ such that $\omega (f)=0$ for all $f$ vanishing at $\infty$, the Dixmier trace of  $T \geq 0$ is defined as 
\[ 
\text{Tr}_{\omega}(T) = \omega (\tau_{\Lambda} (T)).
\]
Then one can extend $\text{Tr}_{\omega}$ to all of $\mathcal{L}^{1,\infty}(\mathcal{H})$ by linearity. 

The resulting linear functional  $\text{Tr}_{\omega}$ is a positive trace on $\mathcal{L}^{(1, \infty)} (\mathcal{H})$ which 
in general depends on the limiting procedure $\omega$. The operators $T \in \mathcal{L}^{1,\infty}(\mathcal{H})$ 
whose Dixmier trace $\text{Tr}_{\omega}(T)$ is independent of the choice of the state $\omega$ are called 
measurable and we will denote their Dixmier trace by $ \int\!\!\!\!\!\!- \, T$. If for a compact positive operator $T$ we have 
\[
\mu_n(T) \sim \frac{c}{n} \qquad {\rm as} \qquad n \to \infty, 
\] 
where $c$ is a constant, then $T$ is measurable and $\int\!\!\!\!\!\!- \,T = c.$
We use this fact to compute the Dixmier trace of  $(1+\triangle')^{-1}$ where 
$\triangle'$ is the perturbed Laplacian on the noncommutative two torus: 

\newtheorem{indixdom}[main1]{Corollary}

\begin{indixdom} \label{indixdom}
The operator $(1+\triangle')^{-1}$ belongs to $\mathcal{L}^{1, \infty}(\mathcal{H}_\varphi)$. 
Moreover it is measurable and  
\[
 \int\!\!\!\!\!\!- \, (1+\triangle')^{-1}= \frac{\pi \varphi(1)}{\Im(\tau)}.
\]

\begin{proof}

Assuming that  the eigenvalues $\lambda_j$ of  $\triangle'$ are written 
in increasing order
\[ 
\lambda_0 \leq \lambda_1 \leq \lambda_2 \leq \cdots,
\]
the asymptotic estimate \eqref{asympcount} can be reformulated as
\[
\lambda_j \sim \frac{\Im(\tau)}{\pi \varphi(1)} j \qquad {\rm as} \qquad  j \to \infty.
\]
Therefore
\[ (1+ \lambda_j)^{-1} \sim  \frac{\pi \varphi(1)}{\Im(\tau)} \big (\frac{\pi \varphi(1)}{\Im(\tau)} + j \big )^{-1} \qquad {\rm as} \qquad  j \to \infty,\]
and it follows that
\[ 
 \int\!\!\!\!\!\!- \,  (1+\triangle')^{-1}  = \frac{\pi \varphi(1)}{\Im(\tau)}.   
\]

\end{proof}

\end{indixdom}

Considering the classical case explained at the beginning of Section \ref{prel}, the 
above theorem suggests the following formula for the volume of the noncommutative two 
torus $\mathbb{T}_\theta^2$ equipped with the metric  associated to the conformal structure 
represented by $\tau$, $\Im(\tau)>0$, and the Weyl factor $k$:
\[ 
\text{Vol}(\mathbb{T}_\theta^2) :=  \frac{4\pi^2}{\Im(\tau)} \varphi(1) = \frac{4\pi^2}{\Im(\tau)}  \mathfrak{t}(k^{-2}).
\]

\section{Connes' Trace Theorem for $\mathbb{T}_\theta^2$} \label{connestracesec} 

Let $M$ be a closed smooth manifold of dimension $n$. Wodzicki defined a trace functional 
on the algebra of pseudodifferential operators of arbitrary order on $M$, and proved that 
it was the only non-trivial trace \cite{wod}. This functional, denoted by $\text{Res}$, is called the 
noncommutative residue. The restriction of $\text{Res}$ to pseudodifferential 
operators of order $-n$ was discovered independently by Guillemin and its properties were studied in 
\cite{guil}. In general, unlike the Dixmier trace, $\text{Res}$ is not a positive linear functional. However 
its restriction to pseudodifferential operators of order $-n$ is positive.  
One of the main results proved in \cite{con0} is that if $E$ is a smooth vector bundle on 
$M$  then the Dixmier trace $\text{Tr}_\omega$ and $\text{Res}$ coincide on pseudodifferential 
operators of order $-n$  acting on $L^2$ sections of  $E$. In fact it is proved that such operators $P$ 
are measurable operators in $\mathcal{L}^{1,\infty}(L^2(M, E))$ and 
\[
 \int\!\!\!\!\!\!- \, P = \frac{1}{n} \text{Res}(P). 
\]
This result is known as Connes' trace theorem.

The noncommutative residue of a classical pseudodifferential operator $P$ acting on smooth sections of 
$E$ is defined as 
\[
\text{Res}(P) = (2 \pi)^{-n} \int_{S^*M} {\rm tr}(\rho_{-n}(x, \xi))\, dx \, d\xi,
\]
where $S^*M \subset T^*M$ is the unit cosphere bundle on $M$ and $\rho_{-n}$ is the component of order 
$-n$ of the complete symbol of $P$. For detailed discussion of traces on classical pseudodifferential operators 
on manifolds and uniqueness results we refer the reader to \cite{lesnei} and references therein.

In \cite{fatwon}, a noncommutative residue on the algebra of  classical pseudodifferential operators on the 
noncommutative two torus is defined and it is proved that up to multiplication by a constant, it is the unique 
continuous trace on the algebra of symbols. A symbol $\rho $ of order $n$ on $A_\theta$ is said to 
be classical if there is an asymptotic expansion of the form 
\[
\rho(\xi) \sim \sum_{j=0}^{\infty} \rho_{n-j}(\xi) \qquad {\rm as} \qquad \xi \to \infty, 
\] 
where  each $\rho_{n-j} : \mathbb{R}^2 \setminus \{0 \} \to A_{\theta}^{\infty}$ is smooth 
and positively homogeneous of order $n-j$. It is shown that the homogeneous terms in 
such an asymptotic expansion are uniquely determined by $\rho$ and the noncommutative residue of 
$P_\rho$ is defined by 
\begin{equation}
{\rm res}(P_\rho) = \int_{\mathbb{S}^1}
{\mathfrak{t}}(\rho_{-2}(\xi))\,d\Omega, \nonumber
\end{equation}
where $d\Omega$ is the Lebesgue measure on the unit circle. This is the analogue of Wodzicki's 
noncommutative residue. 

In the following theorem we prove the analogue of  Connes'  trace theorem \cite{con0} for  
the noncommutative two torus. That is, we show that the Dixmier trace 
and the noncommutative residue defined above coincide on  classical pseudodifferential 
operators of order $-2$ on $\mathbb{T}_\theta^2$. 

\newtheorem{connestrace}[main1]{Theorem}
\begin{connestrace} \label{connestrace}

Let $\rho$ be a classical pseudodifferential symbol of order $-2$ on the noncommutative two torus. 
Then $P_\rho \in \mathcal{L}^{1, \infty}(\mathcal{H}_0)$ and
\[
{\rm Tr}_\omega(P_\rho) = \frac{1}{2}{\rm res}(P_\rho). 
\]

\begin{proof}

We have $P_\rho=A(1+\triangle_0)^{-1}$ where $\triangle_0 = \delta_1^2 + \delta_2^2$ is the Laplacian 
of the metric associated with $\tau=i$ and $k=1$,  and $A=P_\rho(1+\triangle_0)$. Since $A$ is a
pseudodifferential operator of order $0$, it is a bounded operator on $\mathcal{H}_0$. It follows from  
Corollary \ref{indixdom} that $P_\rho$ belongs to the ideal $\mathcal{L}^{1, \infty}(\mathcal{H}_0)$.
Similarly, one can see that any pseudodifferential operator of order $-3$ is a trace-class operator. Since the
Dixmier trace vanishes on trace-class operators, ${\rm Tr}_\omega(P_\rho)$ only depends on the principal part 
of $\rho$. That is, if we write 
\[ \rho(\xi) \sim \rho'(\xi) + \sum_{j=1}^{\infty} \rho_{-2-j}(\xi)  \qquad {\rm as}\qquad \xi \to \infty,\]
where $\rho'$ and $\rho_{-2-j}$ are  positively homogeneous of order $-2$ and $-2-j$,  respectively,  then 
\begin{equation} \label{only-2}
{\rm Tr}_\omega (P_\rho) = {\rm Tr}_\omega (P_{\rho'}).
\end{equation}
Now we prove that 
\begin{equation} \label{onlyitstrace}
{\rm Tr}_\omega(P_{\rho'})={\rm Tr}_\omega(P_{\rho'_{0,0}}),
\end{equation}
where 
\[\rho'_{0,0}:= \mathfrak{t} \circ \rho': \mathbb{R}^2 \setminus \{0\} \to \mathbb{C}.\]
First we make the following observations. The symbol of $P_{\rho'} P_U - P_U P_{\rho'}$ is equivalent to
\[
\rho' U - U\rho' + (\partial_1\rho') U + \frac{1}{2!}  (\partial_1^2\rho') U + \frac{1}{3!} (\partial_1^3\rho') U + \cdots.
\]
Since 
\[
(\partial_1\rho') U + \frac{1}{2!}  (\partial_1^2\rho') U + \frac{1}{3!} (\partial_1^3\rho') U + \cdots
\]
is equivalent to a symbol of order $-3$, the Dixmier trace vanishes on the corresponding operator. 
Therefore, since the Dixmier trace is a trace, we have
\[
{\rm Tr}_\omega (P_{\rho' U} -  P_{U\rho'}) = {\rm Tr}_\omega (P_{\rho'} P_U - P_U P_{\rho'}  ) =0. 
\] 
Hence 
\begin{equation} \label{commwU}
{\rm Tr}_\omega (P_{\rho' U} ) = {\rm Tr}_\omega (P_{ U \rho'}).
\end{equation} 
Using a similar argument for $V$ instead of $U$, we have
\begin{equation} \label{commwV}
{\rm Tr}_\omega (P_{\rho' V} ) = {\rm Tr}_\omega (P_{ V \rho'} ).
\end{equation}
Now, in order to prove \eqref{onlyitstrace}, we consider the expansion of $\rho'(\xi)$ with rapidly 
decaying coefficients:
\[
\rho'(\xi) = \sum_{m, n \in \mathbb{Z}} \rho'_{m, n} (\xi)  U^m V^n.
\] 
Note that since $\rho'$ is positively homogeneous of order -2, so are the functions 
$\rho'_{m,n}: \mathbb{R}^2 \setminus \{0\} 
\to \mathbb{C}$ for all $m, n \in \mathbb{Z}$. Using  \eqref{commwU} and \eqref{commwV} we see that for fixed $m, n \in \mathbb{Z}$, 
the pseudodifferential operator associated with $\rho'_{m, n}   U^m V^n $ has the same Dixmier trace as the ones associated with the following symbols:
\[ 
\rho'_{m, n}   U^{m-1} V^nU= e^{2 \pi i n \theta}\rho'_{m,n} U^{m} V^n, 
\qquad V \rho'_{m, n}   U^{m} V^{n-1} = e^{2 \pi i m \theta} \rho'_{m, n} U^{m} V^n. 
\] 
Since $\theta$ is irrational, it follows that 
\[ 
{\rm Tr}_\omega( P_{\rho'_{m,n} U^m V^n}) = 0 \qquad {\rm if} \qquad (m, n) \neq (0,0). 
\]
Therefore, in order to conclude \eqref{onlyitstrace}, it suffices to show that if  we define 
\[
q(\xi) = \rho'(\xi) - \sum_{|m|, |n| \leq M }   \rho'_{m, n}(\xi)   U^m V^n, 
\]
then ${\rm Tr}_\omega(P_{q})$ goes to $0$ as $M \to \infty$. This can be proved as follows. We write 
$P_q = B (1 + \triangle_0)^{-1}$
where $B$ is a pseudodifferential operator of order $0$, and we denote its symbol by $b$. 
Since $B=P_q (1+ \triangle_0)$ we have 
$
b \sim q (1 + \xi_1^2+\xi_2^2).
$ 
Therefore 
\[
|{\rm Tr}_\omega(P_q)|=|{\rm Tr}_\omega \big (P_{q(1+|\xi|^2)} (1+\triangle_0)^{-1} \big )| 
\leq ||P_{q(1+|\xi|^2)} ||_{\mathcal{L}(\mathcal{H}_0)}  ||(1+\triangle_0)^{-1} ||_{1,\infty}.
\]
Since $q$ is positively homogeneous of order $-2$, as $M$ goes to infinity, the 
$A_\theta$-norm of $q(1+|\xi|^2)$ and its 
derivatives go to $0$ uniformly in $|\xi| \geq 1$. So $||P_{q(1+|\xi|^2)} ||_{\mathcal{L}(\mathcal{H}_0)}$  
approaches $0$ as $M \to \infty$. Hence 
\begin{equation} \label{contwrtsymb}
\lim_{M\to \infty}{\rm Tr}_\omega(P_q) = 0,
\end{equation}
and the identity \eqref{onlyitstrace} is proved.  

Combining \eqref{only-2} and \eqref{onlyitstrace}, so far we have proved that 
\[
{\rm Tr}_\omega(P_\rho)={\rm Tr}_\omega(P_{\rho'})={\rm Tr}_\omega(P_{\rho'_{0,0}}).
\]
Note that $\rho'_{0,0}= \mathfrak{t} \circ \rho'$ is a complex-valued function that is viewed as 
a symbol of order $-2$ on $\mathbb{T}_\theta^2$. This map  is determined by its values on 
$\mathbb{S}^1$ since it is positively homogeneous of order $-2$. 
In order to analyze  ${\rm Tr}_\omega(P_{\rho'_{0,0}})$ we define a linear functional $\mu$ 
on the space of smooth complex-valued functions on $\mathbb{S}^1$ as follows.  Given a 
smooth map $f :\mathbb{S}^1 \to \mathbb{C} \subset A_\theta^\infty$, which will be viewed as an 
$A_\theta$-valued map, first, we extend it to a homogeneous map of order $-2$ from 
$\mathbb{R}^2 \setminus \{ 0\}$ to $A_\theta^\infty$, and denote the extension by $\tilde{f}$. Then 
we define $\mu(f)$ to be the Dixmier trace of the pseudodifferential operator corresponding to $\tilde{f}$: 
\[ 
\mu(f) := {\rm Tr}_\omega(P_{\tilde{f}}), \qquad f \in C^{\infty}(\mathbb{S}^1). 
\]

Considering the positivity of the Dixmier trace and the formulas for product and adjoint 
of pseudodifferential operators explained in Proposition \ref{symbolcalculus}, 
it follows that $\mu$ is a positive linear functional. The argument used for proving identity 
\eqref{contwrtsymb} shows that  $\mu$ extends by continuity to a positive linear functional on 
continuous functions on $\mathbb{S}^1$. Therefore it follows from the Riesz representation theorem 
that $\mu$ is given by integration against a Borel measure on $\mathbb{S}^1$.

Now we show that $\mu$ is rotation invariant. Since the range of any smooth map $f :\mathbb{S}^1 \to \mathbb{C} 
\subset A_\theta^\infty$  is in the center of $A_\theta^\infty$, in the definition of $P_{\tilde{f}}$ given 
by \eqref{pseudofourierdef} the noncommutativity of the algebra does not play a role. So  we can think 
of $P_{\tilde{f}}$ as a pseudodifferential operator on the ordinary two torus. By a straightforward 
computation one can see that if $\sigma(x, \xi)$ is an ordinary pseudodifferential symbol on 
$\mathbb{R}^2$ and if $T$ is a rotation of the plane then
\[
 P_{\sigma(Tx, T_\xi)} = \mathcal{U}_T^{-1} P_{\sigma(x, \xi)}    \mathcal{U}_T, 
\]
where   $\mathcal{U}_T$ is the unitary operator given by
\[ 
\mathcal{U}_T(g) := g \circ T^{-1}, \qquad  g \in C_c^\infty(\mathbb{R}^2).
\] 
This fact and the trace property of the Dixmier trace imply that 
\[
{\rm Tr}_{\omega}(P_{\sigma(Tx, T\xi)})= {\rm Tr}_{\omega}(\mathcal{U}_T^{-1} P_{\sigma(x, \xi)}    \mathcal{U}_T) 
= {\rm Tr}_{\omega}(P_{\sigma(x, \xi)}).  
\]
In case $\sigma(x, \xi)=\tilde{f}(\xi)$, this implies that $\mu(f \circ T) = \mu(f)$. 

Since the linear functional $\mu$ is invariant under rotations of the circle, it 
is given by integration against a constant multiple of the Lebesgue measure on 
$\mathbb{S}^1$. Denoting this constant by $c$, and using \eqref{only-2} and \eqref{onlyitstrace}, 
we can write:
\begin{eqnarray}
{\rm Tr}_\omega (P_\rho) &=& {\rm Tr}_\omega (P_{\rho'})   = {\rm Tr}_\omega (P_{\rho'_{0, 0}}) = 
\mu \Big ({\rho'_{0,0}}_{  {\big \arrowvert}_{\mathbb{S}^1}} \Big ) \nonumber \\
&=&c \int_{\mathbb{S}^1} \rho'_{0,0} \,d\Omega 
=  c \, {\rm res}(P_\rho). \nonumber
\end{eqnarray}
Finally the constant $c$ is fixed by $(1+ \triangle_0)^{-1}$ as follows. Using Corollary 
\ref{indixdom}, we have
$ 
{\rm Tr}_\omega ((1+ \triangle_0)^{-1}) =  \pi.
$
On the other hand one can easily see that
$
{\rm res}((1+ \triangle_0)^{-1}) =2 \pi,
$
because the asymptotic expansion with homogeneous terms for its symbol starts with $1/(\xi_1^2+\xi_2^2)$.
Hence
$
c= \frac{1}{ 2}.
$

\end{proof}

\end{connestrace}

\noindent Department of Mathematics and Statistics, 
York University, 
Toronto, Ontario, Canada, M3J 1P3\\
ffathiza@mathstat.yorku.ca\\

\noindent Department of Mathematics, 
 The University of Western Ontario, 
 London, Ontario, Canada, N6A 5B7\\
 masoud@uwo.ca

\end{document}